\documentclass[11pt]{amsart}
\usepackage{amssymb,amsfonts}

\newtheorem{theorem}{Theorem}[section]
\newtheorem{lemma}[theorem]{Lemma}
\newtheorem{proposition}[theorem]{Proposition}

\newtheorem{definition}[theorem]{Definition}
\newtheorem{corollary}[theorem]{Corollary}

\newcommand{\be}{\begin{equation}}
\newcommand{\ee}{\end{equation}}

\newcommand{\diag}{\mbox{diag}}
\newcommand{\norm}[1]{\Vert #1 \Vert}

\begin{document}

\title{You Can See the Arrows In a Quiver Operator Algebra}

\author{Baruch Solel}

\thanks {partially supported by the Fund for the
 Promotion of Research at the Technion}

\date{}

\begin{abstract} \noindent
We prove that two quiver operator algebras can be isometrically
isomorphic only if the quivers (=directed graphs) are isomorphic.
We also show how the graph can be recovered from certain
representations of the algebra.\\ \textit{2000 subject
classification :} 47L55, 47L75, 46K50.
\end{abstract}

\maketitle

\begin{section}{Introduction}
A quiver is a directed graph with n vertices $\{1,2,\ldots n\}$
and $C_{ij}$ arrows from $j$ to $i$. Here $C_{ij}$ is a non
negative integer (or possibly $\infty$ if the graph is infinite).
Let $A$ be the $C^*$-direct sum of $n$ copies of $\mathbb{C}$
indexed by the vertices. For a finite graph we shall view $A$ as
the algebra $D_n$ of all diagonal $n\times n$ matrices.

As we explain in the next section, one can associate with the
quiver a correspondence $E(C)$ over $A$ and this correspondence
gives rise to a (non selfadjoint) operator algebra, denoted
$\mathcal{T}_+(C)$, that is referred to as the \emph{quiver
algebra}. Another algebra associated with the quiver is
$H^{\infty}(C)$, the $w^*$-closure of $\mathcal{T}_+(C)$. (Here
$C$ is the $n\times n$ matrix $\{C_{ij}\}$ associated with the
quiver). In fact, $\mathcal{T}_+(C)$ is the tensor algebra
associated with the correspondence $E(C)$. (See \cite{tensor} for
more about tensor algebras and their relation to Cuntz-Pimsner
algebras).

If $C$ is the $1\times 1$ matrix whose entry is $n$, we have
$E(C)=\mathbb{C}^n$ and the quiver algebra $\mathcal{T}_+(C)$ is
the non commutative disc algebra $\mathcal{A}_n$ introduced and
studied by Popescu (\cite{popshift} and \cite{An} ). The algebra
$H^{\infty}(C)$ in this case was studied by Popescu
(\cite{popshift}) who denoted it $F^{\infty}_n$ and by Davidson
and Pitts (\cite{DP}) who wrote $\mathcal{L}_n$ for it and referred
 to it as a free semigroup algebra. We shall
use the notation $\mathcal{A}_n$ and $\mathcal{L}_n$ for these
special quiver algebras. For more general (countable) graphs the
algebra that we write as $H^{\infty}(C)$ was recently studied by
Kribs and Power (\cite{KP}). They denoted it $\mathcal{L}_G$ where
$G$ is the graph (=quiver) and called it a free semigroupoid
algebra.

In \cite{An} Popescu proved that, if $n \neq m$ then the algebras
$\mathcal{A}_n$ and $\mathcal{A}_m$ are not isomorphic. One can
also show that, in this case, there is no isomorphism from
$\mathcal{L}_n$ onto $\mathcal{L}_m$ preserving the weak
topologies (that can be deduced from \cite[Theorem 2.3]{DP}). For
the algebras $\mathcal{T}_+(C)$ and $H^{\infty}(C)$ it is shown in
\cite[Corollary 10.5 and Theorem 10.1]{KP} that, if the graphs are
not isomorphic then the algebras cannot be unitarily equivalent .
(With respect to  'regular representations' of the algebras).

The purpose of this note is to relax the condition of
\emph{unitariliy equivalent} and replace it with
\emph{isometrically isomorphic} (for $\mathcal{T}_+(C)$) and
\emph{isometrically isomorphic via a $w^*-w^*$-homeomorphism} (for
$H^{\infty}(C)$). This is proved in Theorem~\ref{main} which is
the main result. Moreover, we show how to 'read' the quiver (that
is, the numbers $C_{ij}$ ) from the characters and certain
representations of the algebra (justifying the title of this
paper).

All the graphs in this paper will be assumed to be countable.
 To simplify the arguments we prove the main result for
graphs with finitely many vertices (that is, $n<\infty$ ) but one
 can extend the arguments (with some care) to general
countable graphs.

The $C^*$-analogue of quiver algebras is referred to as graph
$C^*$-algebras and these have been studied extensively starting
with the work of Cuntz and Krieger \cite{CK} .(See also
\cite{BB,EL,FLR,KPR} and others). It is not true that, if two
graph $C^*$-algebras are isomorphic, then the graphs are
isomorphic. Somehow the non selfadjoint algebra preserves all the
data while the $C^*$-algebra 'forgets' some. A similar phenomenon
was observed by Arveson for algebras associated with dynamical
systems (\cite{Arv}).

\end{section}
\begin{section}{Preliminaries}

  We begin by recalling the notion of a W*-correspondence. For
the general theory of Hilbert C*-modules which we use, we will
follow \cite{Lance}. In particular, a Hilbert C*-module will be a
\emph{right} Hilbert C*-module.
\begin{definition}
\emph{Let $A$ be a von Neumann algebra and let $E$ be a (right)
Hilbert C*-module over $A$. Then $E$ is called a Hilbert W*-module
over $A$ in case it is self dual (that is, every continuous
$A$-module map from $E$ to $A$ is implemented by an element of
$E$). It is called a W*-correspondence over $A$ if it is also
endowed with the structure of a left $A$-module via a normal
*-homomorphism $\varphi : A \rightarrow \mathcal{L}(E)$. (Here
$\mathcal{L}(E)$ is the algebra of all bounded, adjointable,
module maps on $E$)}.
\end{definition}

Given a W*-correspondence over $A$, we denote the $A$-valued inner
product on $E$ by $\langle \cdot,\cdot \rangle$. The \emph{full}
Fock space over $E$ will be denoted by $\mathcal{F}(E)$, so
$$\mathcal{F}(E)=A \oplus E \oplus E^{\otimes 2} \oplus \cdots .$$
(The tensor products here are internal tensor products, see
\cite{Lance}). The space $\mathcal{F}(E)$ is evidently a Hilbert
W*-correspondence over $A$ with left action $\varphi_{\infty}$
given by the formula $$\varphi_{\infty}(a)=\diag (a,\varphi(a),
\varphi^{(2)}(a), \ldots ), $$ where $$\varphi^{(k)}(a)(\xi_1
\otimes \xi_2 \otimes \cdots \xi_k)=\varphi(a)\xi_1 \otimes \xi_2
\otimes \cdots \xi_k .$$ For $\xi \in E$, we write $T_{\xi}$ for
the creation operator on $\mathcal{F}(E)$: $$ T_{\xi}(\xi_1
\otimes \xi_2 \otimes \cdots \xi_k)=\xi \otimes \xi_1 \otimes
\xi_2 \otimes \cdots \xi_k .$$
 Then $T_{\xi}$ is a continuous,
adjointable operator in $\mathcal{L}(\mathcal{F}(E))$. The norm
closed subalgebra of $\mathcal{L}(\mathcal{F}(E))$ generated by
all the $T_{\xi}$'s and $\varphi_{\infty}(A)$ is called \emph{the
tensor algebra} over $E$ and is denoted $\mathcal{T}_+(E)$
(\cite{tensor}). Since $\mathcal{F}(E)$ is a Hilbert W*-module, it
is known that $\mathcal{L}(\mathcal{F}(E))$ is a von Neumann
algebra. We can now close $\mathcal{T}_+(E)$ in the $w^*$-topology
. This $w^*$-closed algebra will be denoted $H^{\infty}(E)$ and
will be referred to as the \emph{weak tensor algebra} of $E$.

We will be interested in a certain class of W*-correspondences.
First, let the algebra $A$ be the algebra $D_n$ of all diagonal
$n\times n$ (complex) matrices (if $n=\infty$ let $A=l_{\infty}$).
For a fixed $n$ let $C$ be a fixed $n\times n$ matrix with entries
in $\mathbb{Z}_+ \cup \{\infty \}$. For each $1\leq i,j \leq n$,
let $E(C)_{ij}$ be a (complex) $C_{ij}$-dimensional Hilbert space.
(We will usually write it as $\mathbb{C}^{C_{ij}}$, where
$\mathbb{C}^{\infty}$ is, of course, $l_2$ and $\mathbb{C}^0$ will
 be understood as $\{0\}$).
The space $E=E(C)$ is the vector space of all $n\times n$ matrices
$\xi$ with the property that its $i,j$ entry, $\xi_{ij}$, is a
vector in $E(C)_{ij}$.( If $n=\infty$, we shall also require that
$$\sup_j \sum_{i=1}^{\infty} \norm{\xi_{ij}}^2 < \infty $$ holds).
 This space can be viewed as an $A-A$
bimodule via the formulae: $$ (\xi D)_{ij}= (\xi)_{ij}d_j $$ $$
(\varphi(D)\xi)_{ij}=d_i (\xi)_{ij}, $$
 where $D=\diag (d_1,d_2,
\ldots )$ lies in $A$. Also, $E(C)$ has an $A$-valued inner
product defined by the formula $$(\langle \xi,\eta \rangle
)_j=\sum_{i=1}^n \langle \xi_{ij},\eta_{ij} \rangle .$$ (Note
that, since the $A$-valued inner product on $E(C)$ is linear in
the \emph{second} term, we shall use this convention also for the
inner products of the Hilbert spaces $E(C)_{ij}$).

This makes $E(C)$ a W*-correspondence over $A$.

Given $E(C)$ as above, we shall write $\mathcal{T}_+(C)$ for
$\mathcal{T}_+(E(C))$ and $H^{\infty}(C)$ for $H^{\infty}(E(C))$.

Note that, when $C$ is the $1\times 1$ matrix whose entry is $n$,
we can write $E(C)=\mathbb{C}^n$. In this case
the algebra $\mathcal{T}_+(C)$ is the algebra $\mathcal{A}_n$
studied by Popescu in \cite{popshift} and in \cite{An} and the algebra
 $H^{\infty}(C)$ is
the algebra $\mathcal{L}_n$ studied by Davidson and Pitts
(\cite{DP}) and by Popescu in \cite{popshift}. (Popescu denoted it
$F_n^{\infty}$). For a general matrix $C$ the algebra
$H^{\infty}(C)$ was studied recently by Kribs and Power
(\cite{KP}). They called it a free semigroupoid algebra and wrote
it $\mathcal{L}_G$ where $G$ is the (countable) graph associated with $C$.
They also studied norm closed algebras $\mathcal{A}_G$
(\cite[Corollary 10.5]{KP}). These are tensor algebras
$\mathcal{T}_+(E(C))$ but, if $n=\infty$, one has to define the
correspondence $E(C)$ as a $C^*$-correspondence over the
$C^*$-algebra $c$ (or $c_0$), not over $l_{\infty}$.

The representation theory for the tensor algebras was worked out
in \cite{tensor}. We now describe some of the basic results.

\begin{definition}
Let $E$ be a W*-correspondence over $A$ and let $H$ be a Hilbert
space.
\begin{itemize}
\item[(1)] A covariant representation of $E$ on $H$ is a pair
$(T,\sigma)$, where
       \begin{itemize}
       \item[(i)] $\sigma$ is a representation of $A$ in $B(H)$.
       \item[(ii)] $T:E \rightarrow B(H)$ is a linear contraction.
       \item[(iii)] $T$ is a bimodule map in the sense that
       $T(\varphi(a)\xi b)=\sigma(a)T(\xi)\sigma(b)$ for
        $\xi \in
       E,; a,b \in A$.
       \end{itemize}
\item[(2)] We say that the covariant representation is completely
contractive if $T$ is.
\item[(3)] We say that the covariant representation is normal if
$\sigma$ is a normal representation and if $T$ is continuous with
respect to the $\sigma$-topology of \cite{BDH} on $E$ and the
$\sigma$-weak topology on $B(H)$.
\item[(4)] We say that the covariant representation is isometric
if, for all $\xi,\eta$ in $E$, $T(\xi)^*T(\eta)=\sigma(\langle
\xi,\eta \rangle )$.
\end{itemize}
\end{definition}

Given a covariant representation $(T,\sigma)$ of $E$ on $H$, we can
define a linear map $\tilde{T}$ from the algebraic tensor product
$E\odot_{\sigma} H$ to $H$ defined by $\tilde{T}(\xi \otimes
h)=T(\xi)h $.

We have the following. (See \cite[Lemma 3.5 and Theorem
3.10]{tensor} for the proof).

\begin{proposition}\label{rho}
\begin{itemize}
\item[(1)] If $\rho$ is a contractive representation of
$\mathcal{T}_+(E)$ on $H$ then, setting
$\sigma(a)=\rho(\varphi_{\infty}(a))$ for $a\in A$ and
$T(\xi)=\rho(T_{\xi})$ for $\xi \in E$, the pair $(T,\sigma)$ is a
covariant representation of $E$.
\item[(2)] The map $\tilde{T}$ defined above is bounded if and
only if $T$ is completely bounded. In fact
$\norm{\tilde{T}}=\norm{T}_{cb}$. So that $\tilde{T}$ is a
contraction if and only if $(T,\sigma)$ is completely contractive.
In this case we view $\tilde{T}$ as a map on the completion
$E\otimes_{\sigma} H$.
\item[(3)] If $(T,\sigma)$ is completely contractive then the
converse to part (1) also holds; i.e. there is a completely
contractive representation $\rho= T\times \sigma $ of
$\mathcal{T}_+(E)$ on $H$ such that
$\sigma(a)=\rho(\varphi_{\infty}(a))$ for $a\in A$ and
$T(\xi)=\rho(T_{\xi})$ for $\xi \in E$.
\end{itemize}
\end{proposition}

In addition to the map $\tilde{T}$ we also define the maps
$\tilde{T}_k :E^{\otimes k}\otimes H \rightarrow H$ by
$\tilde{T}_k(\xi_1 \otimes \cdots \otimes \xi_k)=T(\xi_1) \cdots
T(\xi_k)h$ and then we have
$\tilde{T}_{k+1}=\tilde{T}(I_E\otimes \tilde{T}_k) $.

Given a representation $\pi_0$ of $A$ on a Hilbert space $H_0$ we can
form the Hilbert space $\mathcal{F}(E) \otimes_{\pi_0} H_0$ (where
the inner product is given by $\langle X\otimes h,Y \otimes g
\rangle = \langle h,\pi_0(\langle X,Y \rangle) g\rangle $ for $X,Y$ in
$\mathcal{F}(E)$ and $h,g \in H_0$) and define an isometric
covariant representation $(V,\pi)$ of $E$ on this Hilbert space by
$V(\xi)=T_{\xi} \otimes I_{H_0}$ and $\pi(a)=\varphi_{\infty}(a)
\otimes I_{H_0} $. Such a representation is said to be
\emph{induced}.

The resulting representation $V\times \pi$ of $\mathcal{T}_+(E)$ is, in
fact, the restriction to $\mathcal{T}_+(E)$ of the representation
induced by $\pi_0$ of $\mathcal{L}(\mathcal{F}(E))$ on this
Hilbert space given by the formula
$$ \pi_0^{\mathcal{F}(E)}(T)=T\otimes I_{H_0} \;,\;\;\;\;\;
 T\in \mathcal{L}(\mathcal{F}(E)).$$
 This shows that,
when $(V,\pi)$ is a normal induced representation of $E$,
 the representation $V\times \pi$ can be
extended to a $w^*$-continuous representation of $H^{\infty}(E)$.
In \cite[Proposition 2.8]{curvature} it is shown that, if
$(T,\sigma)$ is a normal completely contractive representation of
$E$ satisfying $\tilde{T}_k\tilde{T}_k^* \rightarrow 0 $ in the
strong operator topology, then the minimal isometric dilation $(V,\pi)$ of
$(T,\sigma)$ is an induced representation (and, thus, the
associated representation $V \times \pi$ of $\mathcal{T}_+(E)$ can be extended to
a $w^*$-continuous representation of $H^{\infty}(E)$).
Since $T\times \sigma$ is then a compression of $V \times \pi$, if
follows that we can also extend this representation to a
$w^*$-continuous representation of $H^{\infty}(E)$. We summarize
this discussion as follows.

\begin{lemma}\label{pure}
If $(T,\sigma)$ is a normal completely contractive covariant
representation of $E$ such that $\tilde{T}_k\tilde{T}^*_k
\rightarrow 0 $ in the strong operator topology then the
representation $T\times \sigma$ can be extended to a
$w^*$-continuous representation of $H^{\infty}(E)$.
\end{lemma}

Restricting to the case $E=\mathbb{C}^n$ (and
$H^{\infty}(E)=\mathcal{L}_n$), we have the following.

\begin{lemma}\label{shift}
Suppose $V=(V_1,V_2, \ldots ,V_n)$ is an n-tuple of isometries in
$B(H)$ (where we allow $n=\infty$) whose ranges are orthogonal and
the sum of the ranges is not all of $H$. Let $\rho$ be the
representation of $\mathcal{T}_+(\mathbb{C}^n)$ ($=\mathcal{A}_n$)
defined by $V$. Then the following hold:
\begin{itemize}
\item[(1)] There is a Hilbert space $K$ and a unitary operator $v:
H \rightarrow \mathcal{F}(\mathbb{C}^n) \otimes K $ such that $$
vV_iv^*=T_{e_i} \otimes I_K
\;(=\pi_0^{\mathcal{F}(\mathbb{C}^n)}(T_{e_i}) ), $$
where
$\{e_i\}$ is the standard orthonormal basis of $\mathbb{C}^n$ and
$\pi_0$ is the obvious representation of $A=\mathbb{C}$ on $K$.
\item[(2)] $v\rho(\cdot)v^*$ is the restriction of
$\pi_0^{\mathcal{F}(\mathbb{C}^n)}$ to
$\mathcal{T}_+(\mathbb{C}^n) \;(=\mathcal{A}_n) $.
\item[(3)] $\rho$ can be extended to a completely isometric
isomorphism of $\mathcal{L}_n$ into $B(H)$ that is a
$w^*-w^*$-homeomorphism onto its image.
\end{itemize}
\end{lemma}
\begin{proof}
In \cite{popshift}, n-tuples as above were called orthogonal
shifts and part (1) follows from Theorem 1.2 there. One can also
deduce it from the above discussion since $(V,\pi)$ is an induced
representation of $E=\mathbb{C}^n$ (where $\pi$ is the
 obvious representation of $\mathbb{C}$ on $H$).
Part (2) follows immediately from (1) and for part (3) note that
the representation $\pi_0^{\mathcal{F}(\mathbb{C}^n)}$ is the
representation that maps $S \in B(\mathcal{F}(\mathbb{C}^n)) $ to
$S\otimes I_K$ in $B(\mathcal{F}(\mathbb{C}^n) \otimes K)$. Since
this representation is completely isometric and a
$w^*-w^*$-homeomorphism of the von Neumann algebra
 $B(\mathcal{F}(\mathbb{C}^n))$ onto its image, the same holds for the
restriction to $\mathcal{L}_n$.
\end{proof}

\end{section}
\begin{section}{Isomorphic Quiver Algebras}
     In this section we prove the main results. We now fix $n<\infty$
 and an $n\times n$ matrix $C$ with entries in $\mathbb{Z}_+\cup \{\infty\}$.
Let $A$ be the C*-algebra $D_n$ of all diagonal $n\times n$
matrices
and let $E(C)$ be the $A-A$-correspondence defined by $C$.

We start by identifying the characters of $\mathcal{T}_+(C)$ and of
$H^{\infty}(C)$. In order to do it we shall first embed $\mathcal{A}_{C_{ii}}$ in
$\mathcal{T}_+(C)$ and $\mathcal{L}_{C_{ii}}$ in $H^{\infty}(C)$.

 Write
$\pi$ for the usual representation of $A$ on $\mathbb{C}^n$
 and write $H$ for $\mathbb{C}^n$.
Let $K(C)$ denote the representation space of the representation
$\pi^{\mathcal{F}(E(C))}$, induced by $\pi$, that is, $$K(C)=
\mathcal{F}(E(C)) \otimes_{\pi} H .$$ For every $1\leq i \leq n$
with $C_{ii}\neq 0$, write $P_i$ for the projection
$\varphi_{\infty}(e_i)\otimes I_H$ in $B(K(C))$ and $K_i(C)$ for
its range. Hence
$$K_i(C)=\varphi_{\infty}(e_i)(\sum_{k=0}^{\infty} \oplus
E(C)^{\otimes k} \otimes H) .$$ For every $1\leq i \leq n$, let
$\{e_{ii}^{(j)} : 1 \leq j \leq C_{ii} \;\}$
 be an orthonormal
 basis for $E(C)_{ii}$ and view these vectors as elements of $E(C)$.
  For $1 \leq j \leq C_{ii}$ let $V_j$ be the operator on
$K_i(C)$ defined by
$$ V_j= \pi^{\mathcal{F}(E(C))}(e_{ii}^{(j)}) | K_i(C) .$$
(Note that $\pi^{\mathcal{F}(E(C))}(e_{ii}^{(j)})$ vanishes on the orthogonal
complement of $K_i(C)$.)
We get $C_{ii}$ isometries on $K_i(C)$ satisfying the conditions
 of Lemma~\ref{shift}. Letting $\Psi_i$ be the map $\rho$ of that
  lemma (with $H=K_i(C)$) composed with the embedding of
  $B(K_i(C))$ into $B(K(C))$ (by defining the operator to be zero
  on the orthogonal complement of $K_i(C)$), we get the following.
  (Note that $H^{\infty}(C)$ can be identified with its image
  under $\pi^{\mathcal{F}(E(C))}$).

\begin{proposition}\label{Psi}
For every $1\leq i \leq n$ with $C_{ii} \neq 0$ there is a
 (non unital) completely isometric isomorphism
$\Psi_i$ of $\mathcal{L}_{C_{ii}}$ into $H^{\infty}(C)$ that is a
$w^*$-homeomorphism (onto its image) and that restricts to a
completely isometric isomorphism of $\mathcal{A}_{C_{ii}}$ into
$\mathcal{T}_+(E(C))$ (denoted also $\Psi_i$).
\end{proposition}

A character of $\mathcal{T}_+(E(C))$ is a one dimensional
(completely) contractive representation and , as such, it is given
by a completely contractive covariant representation $(T,\tau )$
of the W*-correspondence $E(C)$. Here $\tau$ is a one dimensional
representation of $A=D_n$ and, thus, is $\delta_i$ for some $1\leq
i \leq n$ (where $\delta_i$ of a diagonal matrix
$D=diag(d_1,d_2,\ldots ,d_n)$ is $d_i$). The map $T$ is a
contractive linear functional on $E(C)$ satisfying $T(D_1 \xi
D_2)=\delta_i(D_1)T(\xi) \delta_i(D_2) $. Hence it is in fact a
contractive linear functional on the $C_{ii}$-dimensional Hilbert
space $E(C)_{ii}$. If $C_{ii}=0$ then $T=0$. Otherwise,
identifying
 $E(C)_{ii}$ with
$\mathbb{C}^{C_{ii}}$, we associate with every such character a
pair $(i,\lambda)$ where $1\leq i \leq n$ and $\lambda$ is in the
closed unit ball $\overline{\mathbb{B}_{C_{ii}}}$ of
$\mathbb{C}^{C_{ii}}$. If this character can be extended to a
$w^*$-continuous character of $H^{\infty}(C)$ then, using the map
$\Psi_i$ of Proposition~\ref{Psi}, it induces a $w^*$-continuous
character on $\mathcal{L}_{C_{ii}}$. It follows from Theorem 2.3
of \cite{DP} that $\lambda$ lies in the \emph{open} unit ball. We
summarize the discussion in the following theorem. To simplify the
statement, we shall assume that, whenever $C_{ii}=0$, the notation
$\mathbb{B}_{C_{ii}}$ (or its closure) will be interpreted as
$\{0\}$. When $C_{ii}=\infty$ the balls $\mathbb{B}_{C_{ii}}$ and
$\overline{\mathbb{B}_{C_{ii}}}$ are the balls of $l^2$ equipped
with the weak topology.

\begin{theorem}\label{char}
Let $M(\mathcal{T}_+(C))$ be the set of all characters of the
quiver algebra $\mathcal{T}_+(C)$ equipped with the $w^*$-topology
and let $M(H^{\infty}(C))$ be the set of all $w^*$-continuous
characters of $H^{\infty}(C)$ (also with the $w^*$-topology).
Then we have the following homeomorphisms:
\begin{itemize}
\item[(1)] $M(\mathcal{T}_+(C)) \cong \bigcup \{(i,\lambda) :
\lambda \in \overline{\mathbb{B}_{C_{ii}}}, 1\leq i \leq n\} .$
\item[(2)] $M(H^{\infty}(C)) \cong \bigcup \{(i,\lambda) : \lambda
\in \mathbb{B}_{C_{ii}}, 1\leq i \leq n \} .$
\end{itemize}
(Each set on the right hand side is a disjoint union of $n$ closed
and open sets). The character, $\varphi_{(i,\lambda)}$,
 associated with $(i,\lambda)$ is equal $\delta_i$ on $A$ and on
 $T_{\xi}$, for $\xi \in E(C)$, it is defined by
$$\varphi_{(i,\lambda)}(T_{\xi})= \langle \lambda, \xi_{ii}
\rangle .$$
\end{theorem}
\begin{proof} The identifications in both (1) and (2) were shown
above. The fact that these are homeomorphisms is easy to check.
\end{proof}

\vspace{4 mm}

In the following we shall also interpret $\mathbb{C}^m$ and
$\mathbb{B}_m$ as $\{0\}$ if $m=0$.
Now fix $1\leq i,j \leq n$ , $i\neq j$ and $\tilde{\lambda} = (\lambda_i,\lambda_j)
 \in \mathbb{B}_{C_{ii}} \times \mathbb{B}_{C_{jj}} $ and define
  $G(C,\tilde{\lambda},i,j)$
 (respectively, $G_0(C,\tilde{\lambda},i,j)$) to be the set of
all contractive representations $\rho$ of $\mathcal{T}_+(C)$
 (respectively, all contractive $w^*$-continuous
 representations of $H^{\infty}(C)$) on the space $\mathbb{C}^2$
  satisfying:
 \begin{itemize}
 \item[(G1)] For $D \in A$,
 $\rho(D)=\diag (\delta_i(D),\delta_j(D))$.
 \item[(G2)] The image of $\rho$ is contained in $T_2$, the upper
 triangular $2\times 2$ matrices.
 \item[(G3)] For every $S$ in the algebra, $(\rho(S))_{11}=
  \varphi_{(i,\lambda_i)}(S)$ and $(\rho(S))_{22}= \varphi_{(j,\lambda_j)}(S)$.
\end{itemize}

We now present examples of representations in $ G(C,\tilde{\lambda},i,j)$.
 Write $\sigma$ for the representation of $A$ on
$\mathbb{C}^2$ given by $\sigma(D)=\diag
(\delta_i(D),\delta_j(D))$. For $\gamma$ in $\mathbb{C}^{C_{ij}}$,
define the map $T_{\gamma}: E(C) \rightarrow T_2$ by $$
T_{\gamma}(\xi) = \left( \begin{array}{cc}
\varphi_{(i,\lambda_i)}(\xi_{ii}) & \langle \gamma, \xi_{ij}
\rangle \\ 0 & \varphi_{(j,\lambda_j)}(\xi_{jj}) \\ \end{array}
\right)= \left( \begin{array}{cc} \langle \lambda_i, \xi_{ii}
\rangle & \langle \gamma, \xi_{ij} \rangle \\ 0 & \langle
\lambda_j, \xi_{jj} \rangle \end{array} \right)
 \; \in T_2 .$$
Then the pair $(T_{\gamma},\sigma)$ satisfies $T_{\gamma}(D_1\xi
D_2)= \sigma(D_1)T_{\gamma}(\xi)\sigma(D_2) $ for all $D_1,D_2$ in
$A$ and $\xi$ in $E(C)$.

Before we proceed to show that this construction gives us a
representation in $G(C,\tilde{\lambda},i,j)$ we observe that the
converse holds.
\begin{lemma}\label{converse}
Let $\rho$ be a representation in $G(C,\tilde{\lambda},i,j)$ and
let $T(\xi)$ be $\rho(T_{\xi})$ (in $T_2$). Then there is some
$\gamma$ in $\mathbb{C}^{C_{ij}}$ with $\norm{\gamma }^2+
\norm{\lambda_i}^2\leq 1$ such that $T=T_{\gamma}$.
\end{lemma}
\begin{proof}
Since $\rho \in G(C,\tilde{\lambda},i,j)$, it follows that
$T(\xi)_{11}=\varphi_{(i,\lambda_i)}(T_{\xi})=\langle
\lambda_i, \xi_{ii} \rangle$, $T(\xi)_{22}=\langle
\lambda_j, \xi_{jj} \rangle $ and, for $D \in A$,
$\rho(\varphi_{\infty}(D))=\sigma(D)$. Hence $T(D_1\xi
D_2)=\sigma(D_1)T(\xi)\sigma(D_2) $ (for $D_i$ in $A$) and it
follows that $T(\xi)_{12}$ depends only on $\xi_{ij}$. Since this
dependence is clearly linear and $T$ is bounded, we find that
$T(\xi)_{12}$ is $\langle \gamma, \xi_{ij} \rangle $ for some
$\gamma \in \mathbb{C}^{C_{ij}} $. Thus $T=T_{\gamma}$.
To show that $\norm{\gamma}^2+\norm{\lambda_i}^2\leq 1$ let $\xi
\in E(C)$ be defined by : $\xi_{ii}=\lambda_i/\norm{\lambda_i}$,
$\xi_{ij}=\gamma / \norm{\gamma}$ and $\xi_{lp}=0$ otherwise. (The
case where either $\gamma$ or $\lambda_i$ is zero can be handled
easily).
Then $\langle \xi, \xi \rangle $ is the diagonal matrix whose
$i$th diagonal entry is $\norm{\xi_{ii}}^2=1 $, the $j$th diagonal
entry is $\norm{\xi_{ij}}^2=1 $ and all other entries equal zero.
Hence $\norm{T_{\xi}}=\norm{\xi}=1 $ and, consequently,
$\norm{T_{\gamma}(\xi)}=\norm{\rho(T_{\xi})}\leq 1$. But
$$T_{\gamma}(\xi)= \left( \begin{array}{cc} \langle
\lambda_i,\xi_{ii} \rangle & \langle \gamma, \xi_{ij} \rangle \\
0 &0 \end{array} \right) = \left( \begin{array}{cc}
\norm{\lambda_i} & \norm{\gamma} \\0&0 \end{array} \right). $$
Thus $1 \geq \norm{T_{\gamma}(\xi)}^2= \norm{\lambda_i}^2+
\norm{\gamma}^2 $.
\end{proof}

Now recall that the pair $(T_{\gamma},\sigma)$ defines a map
$\tilde{T}_{\gamma}$ from the algebraic tensor product $E(C)
\odot_{\sigma} \mathbb{C}^2$ to $\mathbb{C}^2$ satisfying
$$\tilde{T}_{\gamma} (\xi \otimes h) = T_{\gamma}(\xi)h $$
which is bounded (respectively, contractive) if and only if
 $T_{\gamma}$ is completely bounded (respectively, completely
 contractive). If $\tilde{T}_{\gamma}$ is contractive, then
  (using Proposition~\ref{rho}) the
 pair $(T_{\gamma},\sigma)$ defines a completely contractive
  representation $T_{\gamma}\times \sigma $ of $\mathcal{T}_+(C)$.

 The proof of the following lemma is a straightforward computation
 and is omitted.

\begin{lemma}\label{gamma}
Let $T_{\gamma}$ be as above (for some $\gamma \in
\mathbb{C}^{C_{ij}}$)
 and fix $\xi$ in $E(C)$. Then, for $k=(k_1,k_2)^t $
and $h=(h_1,h_2)^t$ in $\mathbb{C}^2$ we have
\begin{itemize}
\item[(1)] $$ \tilde{T}_{\gamma}(\xi \otimes h)= \left(
\begin{array}{c} \langle \lambda_i, \xi_{ii}\rangle h_1 + \langle
\gamma, \xi_{ij}\rangle h_2 \\ \langle \lambda_j,\xi_{jj}\rangle
h_2 \end{array} \right) .$$
\item[(2)] $$ \tilde{T}_{\gamma}^* k = \eta \otimes \left(
\begin{array}{c} 1 \\1 \end{array} \right) $$
where $\eta$ is the element of $E(C)$ with
$\eta_{ii}=k_1\lambda_i$, $\eta_{ij}=k_1\gamma$,
$\eta_{jj}=k_2\lambda_j$ and all other entries are $0$.
\item[(3)] $$\tilde{T}_{\gamma}\tilde{T}^*_{\gamma}= \left(
\begin{array}{cc} \norm{\lambda_i}^2 + \norm{\gamma}^2 & 0 \\
0 & \norm{\lambda_j}^2  \end{array} \right). $$
\end{itemize}
\end{lemma}

In part (2) of the following corollary we present the general form
of the representations in $G(C,\tilde{\lambda},i,j)$.

\begin{corollary}\label{w*rep}
\begin{itemize}
\item[(1)] $\norm{\tilde{T}_{\gamma}} \leq 1$ (that is, it defines a
completely contractive representation $\rho_{\gamma}=T_{\gamma}\times \sigma$
 of the algebra
$\mathcal{T}_+(C)$ ) if and only if $\norm{\gamma}^2 \leq 1-
\norm{\lambda_i}^2 $.
\item[(2)] The representations in $G(C,\tilde{\lambda},i,j)$ are
all completely contractive and of the form
$\rho_{\gamma}=T_{\gamma}\times \sigma $ (for some $\gamma$ in
$\mathbb{C}^{C_{ij}}$ with $\norm{\gamma}^2 \leq
1-\norm{\lambda_i}^2 $).
\item[(3)] When $\norm{\gamma}^2 \leq 1- \norm{\lambda_i}^2 $, we
have $\norm{\tilde{T}_{\gamma,k}} \rightarrow 0 $. Hence
$\tilde{T}_{\gamma}$ defines a $w^*$-continuous representation of
$H^{\infty}(C)$. Therefore, $G(C,\tilde{\lambda},i,j) =
G_0(C,\tilde{\lambda},i,j) $.
\end{itemize}
\end{corollary}
\begin{proof}
Parts (1) and (2) follow from Lemma~\ref{gamma} and Lemma~\ref{converse}.
 For part (3), fix
$\gamma$ with $\norm{\gamma}^2 \leq 1- \norm{\lambda_i}^2$ and
write $T$ for $T_{\gamma}$. Recall that $\tilde{T}_k$ is a map
from $E^{\otimes k} \otimes \mathbb{C}^2 $ to $\mathbb{C}^2$
defined recursively by $\tilde{T}_1=\tilde{T}$ and
$\tilde{T}_{k+1}=\tilde{T}(I_E \otimes \tilde{T}_k)$. Hence
$\tilde{T}_{k+1}\tilde{T}_{k+1}^* =\tilde{T} (I\otimes
\tilde{T}\tilde{T}^*) \tilde{T}^* $. Thus, if
$$ \tilde{T}_k\tilde{T}_k^* = \left( \begin{array}{cc} a & 0 \\ 0 & b
\end{array} \right) $$
then, for $g=(g_1,g_2)^t$ in $\mathbb{C}^2 $, $$
\tilde{T}_{k+1}\tilde{T}_{k+1}^*g= \tilde{T} \left(I \otimes
\left(
\begin{array}{cc} a&0 \\ 0&b  \end{array} \right) \right)
\tilde{T}^*g= \tilde{T}\left( \eta \otimes \left( \begin{array}{c}
a\\ b
\end{array} \right) \right) $$
where $\eta_{ii}=g_1\lambda_i$, $\eta_{ij}=g_2\gamma$,
$\eta_{jj}=g_2\lambda_j$ and $\eta_{lp}=0$ otherwise.
Hence
$$\tilde{T}_{k+1}\tilde{T}_{k+1}^*g=
\left( \begin{array}{c} \langle \lambda_i,g_1\lambda_i \rangle a
+ \langle \gamma,g_1 \gamma \rangle b \\ \langle \lambda_j,g_2
\lambda_j \rangle b  \end{array} \right) = \left( \begin{array}{c}
g_1a \norm{\lambda_i}^2 + g_1b \norm{\gamma}^2 \\ g_2b
\norm{\lambda_j}^2  \end{array} \right) =$$
$$= \left( \begin{array}{cc}
a\norm{\lambda_i}^2 +b \norm{\gamma}^2 & 0 \\ 0 &
b\norm{\lambda_j}^2 \end{array} \right) g $$
and $$ \tilde{T}_{k+1}\tilde{T}^*_{k+1} = \left( \begin{array}{cc}
a\norm{\lambda_i}^2 + b\norm{\gamma}^2 & 0 \\ 0 &
b\norm{\lambda_j}^2  \end{array} \right) .$$
Write $q_i = \norm{\lambda_i}^2$ and $t=\norm{\gamma}^2$. Then the
 computation above shows
that $\norm{\tilde{T}}^2= \mbox{max}\{q_1 +t, q_2\}$,
$\norm{\tilde{T}_2}^2=\mbox{max}\{q_1^2+tq_1+tq_2, q_2^2 \}$ and,
in general, $$ \norm{\tilde{T}_k}^2=max\{
q_1^k+t(q_1^{k-1}+q_1^{k-2}q_2+ \cdots +q_2^{k-1}),q_2^k \} .$$ If
$q=\mbox{max}\{ q_1,q_2\} $ then $q<1$ and $$ \norm{\tilde{T}_k}
\leq q^k+ktq^{k-1} \rightarrow 0 .$$ Using Lemma~\ref{pure} we see
that this implies that the representation $\rho_{\gamma}$ can be
extended to a $w^*$-continuous representation of $H^{\infty}(C)$.
\end{proof}

We now equip the set $G(C,\tilde{\lambda},i,j)$ with the topology
of pointwise convergence. It is then homeomorphic to a closed
subset of the product space $$ \prod \{\overline{\mathbb{B}}(T_2)
: S\in \mathcal{T}_+(C), \norm{S} \leq 1 \} $$ (equipped with the
product topology), where $\overline{\mathbb{B}}(T_2)$ is the
closed unit ball of $T_2$. This shows that, in this topology,
$G(C,\tilde{\lambda},i,j)$ is a compact set.

For every $\eta$ in $E(C)_{ij}$ ($\cong \mathbb{C}^{C_{ij}}$ ), $$
\langle \gamma, \eta \rangle =(T_{\gamma}\left(
\begin{array}{cc} 0 & \eta \\ 0 & 0 \end{array} \right)  )_{12} =(\rho_{\gamma}
( T_{\tilde{\eta}} ))_{12} $$ where $\tilde{\eta}=\left(
\begin{array}{cc} 0 & \eta \\ 0 & 0 \end{array} \right)$.
 If $\rho_{\gamma_{\alpha}} \rightarrow \rho_{\gamma} $ in
$G(C,\tilde{\lambda},i,j) $ then, for every $\eta \in E(C)_{ij}$,
$\langle \gamma_{\alpha}, \eta \rangle \rightarrow \langle \gamma,
\eta \rangle$, that is, $\gamma_{\alpha} \rightarrow \gamma $ in
the weak topology (of $E(C)_{ij}$). Since the map $\rho_{\gamma}
\mapsto \gamma$ is a bijection and the spaces are compact, we
conclude
 that it is a homeomorphism. We summarize it in the following
proposition.

\begin{proposition}\label{homeo}
The set $G(C,\tilde{\lambda},i.j)$ , equipped with the topology of
pointwise convergence, is homeomorphic to a closed ball (of positive radius) in
$\mathbb{C}^{C_{ij}}$ (equipped with the weak topology).
\end{proposition}

\begin{theorem}\label{main}
Let $C$ be in $M_n(\mathbb{Z}_+\cup \{\infty \})$ and $C'$ be in
$M_m(\mathbb{Z}_+\cup \{\infty\})$.
\begin{itemize}
\item[(1)] If the algebras $\mathcal{T}_+(C)$ and
$\mathcal{T}_+(C')$ are isometrically isomorphic then $n=m$ and
there is a permutation $\tau \in S_n$ such that
$$ C'_{ij}=C_{\tau(i) \tau(j)} $$ for all $i,j$.
\item[(2)] If the algebras $H^{\infty}(C)$ and $H^{\infty}(C')$
are isometrically isomorphic via an isomorphism that is
$w^*$-bicontinuous then $n=m$ and there is a permutation $\tau \in
S_n$ such that
$$ C'_{ij}=C_{\tau(i) \tau(j)} $$ for all $i,j$.
\end{itemize}
\end{theorem}
\begin{proof}
We start by proving (1). Write $\Lambda : \mathcal{T}_+(C)
\rightarrow \mathcal{T}_+(C') $ for the isometric isomorphism. The
selfadjoint part of $\mathcal{T}_+(C)$, $\mathcal{T}_+(C) \cap
\mathcal{T}_+(C)^*$, is equal to $\varphi_{\infty}(A)$ which is
isomorphic to $D_n$. Since $\Lambda$ is an isometry, it maps the
selfadjoint part of $\mathcal{T}_+(C)$ onto the selfadjoint part
of $\mathcal{T}_+(C')$ (\cite{AS}), thus inducing an isomorphism
of $D_n$ onto $D_m$. This shows that $n=m$ and there is some
permutation $\tau \in S_n$ such that, whenever $D$ is in $A$, $$
\Lambda(\varphi_{\infty}(D))=\varphi_{\infty}(\tau^{(n)}(D)) $$
where $\tau^{(n)}(\diag (d_1,\ldots,d_n))=\diag
(\tau(d_1),\ldots,\tau(d_n)) $. Fix a character,
$\varphi=\varphi_{(i,\lambda)} $ in $M(\mathcal{T}_+(C'))$. Then
$\varphi \circ \Lambda $ is in $M(\mathcal{T}_+(C))$. Restricted
to $A$, the character $\varphi =\varphi_{(i,\lambda)}$ vanishes on
diagonal matrices whose $i$th entry is $0$. Thus $\varphi \circ
\Lambda$ vanishes on these diagonal matrices whose
$\tau^{-1}(i)$th entry is zero. We can write $$
\varphi_{(i,\lambda)} \circ \Lambda =
\varphi_{(\tau^{-1}(i),\lambda')} .$$ In fact, it is clear from
Theorem~\ref{char} (part (1)) that $\lambda \mapsto \lambda'$ is a
homeomorphism, denoted $\theta_i$, from
$\overline{\mathbb{B}_{C'_{ii}}}$ onto
$\overline{\mathbb{B}_{C_{\tau^{-1}(i) \tau^{-1}(i)}}} $. Hence,
for every $1\leq i\leq n$, $$ C'_{\tau(i) \tau(i)}=C_{ii} .$$ Now
fix $i \neq j$. For a representation $\rho \in
G(C',\mathbf{0},i,j) $ (where $\mathbf{0}=(0,0) $), $\rho \circ
\Lambda$ is a contractive representation of $\mathcal{T}_+(C)$
whose image is contained in $T_2$. For $D=\diag (d_1,d_2,\ldots
d_n)$ in $A$, $$\rho \circ \Lambda(D)=\rho(\diag
 (\tau(d_1),\ldots
\tau(d_n))) = \diag (\delta_{\tau(i)}(D),\delta_{\tau(j)}(D)) .$$
Also, for $S \in \mathcal{T}_+(C)$, $$\rho \circ \Lambda
(S)_{11}=\rho(\Lambda(S))_{11} = \varphi_{(i,0)}(\Lambda(S))=
\varphi_{(\tau^{-1}(i),\theta_i(0))}(S) $$ and, similarly, $$ \rho
\circ \Lambda(S)_{22}=\varphi_{(\tau^{-1}(j),\theta_j(0))}(S) .$$
Since $\theta_i$ and $\theta_j$ are homeomorphisms of the closed
unit balls and $0$ is an interior point, $\theta_i(0)$ and
$\theta_j(0)$ are in the open unit balls. We write
$\tilde{\lambda}=(\theta_i(0),\theta_j(0))$ and conclude that
$\rho \circ \Lambda $ lies in
$G(C,\tilde{\lambda},\tau^{-1}(i),\tau^{-1}(j)) $.

Since the map $\rho \mapsto \rho \circ \Lambda$ is a homeomorphism
(with respect to the topology of pointwise convergence) we get (using
 Proposition~\ref{homeo}) a
homeomorphism of the closed unit ball in
$\mathbb{C}^{C_{\tau^{-1}(i),\tau^{-1}(j)}}$ and the closed unit ball in
$\mathbb{C}^{C'_{ij}}$. Thus
$$ C_{\tau^{-1}(i),\tau^{-1}(j)} = C'_{ij} .$$
This proves part (1). The proof of part (2) is almost identical
except that in the argument showing $C_{ii}=C'_{\tau(i)\tau(i)}$ we
 use open balls (and part (2) of
Theorem~\ref{char}) instead of closed balls. The proof for $i\neq j$
 is the same due
to the fact that every representation in
$G(C,\tilde{\lambda},i,j)$ extends to a $w^*$-continuous
representation of $H^{\infty}(C)$ ( Corollary~\ref{w*rep} (3)).

\end{proof}
\end{section}

\vspace{3mm}

\noindent
 \small Department of Mathematics \\ Technion
 \\Haifa 32000\\
 Israel\\\texttt{mabaruch@techunix.technion.ac.il}

\end{document}